\documentstyle[12pt]{article}

\begin{document}
\newcommand{\ol }{\overline}
\newcommand{\ul }{\underline }
\newcommand{\ra }{\rightarrow }
\newcommand{\lra }{\longrightarrow }
\newcommand{\ga }{\gamma }
\newcommand{\st }{\stackrel }

\title{On ${\cal N}_c$-Covering Groups of a Nilpotent Product of Cyclic Groups}
\author{by\\ B. Mashayekhy and A. Khaksar\\ Department of Mathematics,\\ Ferdowsi University of
Mashhad,\\ P.O.Box 1159-91775, Mashhad, Iran\\
E-mail: mashaf@math.um.ac.ir }
\date{}
\maketitle
\begin{abstract}
 The first author, in 2001, presented a structure for the
 Baer-invariant of a nilpotent product of cyclic groups with
 respect to the variety of nilpotent groups. In this paper, using
 the above structure, we prove the existence of an ${\cal
 N}_c$-covering group for a nilpotent product of a family of
 cyclic groups. We also present a structure for the ${\cal
 N}_c$-covering group.
\footnote{
A.M.S. Classification: 20F18, 20E34, 20E10.}
\end{abstract}
\newpage
 To find an introduction and preliminaries, we refer the reader to
 see [2,3] and [4]. The following notations are used throughout
 the paper.

  Let $\{A_i|i\in I\}$ be a family of groups where $I$ considered
  as an ordered set. For each $i\in I$, $L_i$ denotes a fixed ${
  \cal N}_c$-covering group for $A_i$ (if any) such that the
  sequence $1\ra M_i\ra L_i \st {\nu_i}{\ra} A_i\ra 1$ is exact
  and
$$ M_i\subseteq Z_c(L_i)\cap \ga_{c+1}(L_i)\ \ and\ \ M_i\cong
{\cal N}_cM(A_i),$$ where $M_i$ is a normal subgroup of $L_i$,
$Z_c(L_i)$ and $\ga_{c+1}(L_i)$ are the $c$th terms of the upper
and the lower central series of $L_i$, respectively, and ${\cal
N}_cM(A_i)$ is the Baer-invariant of $A_i$ with respect to the
variety of nilpotent groups of class at most $c,\ {\cal N}_c$.
Assume that $A=\prod_{i\in I}^{*}A_i$ and $L=\prod_{i\in
I}^{*}L_i$ are free products of $A_i$'s and $L_i$'s, respectively.
We denote by $\nu$ the natural homomorphism from $L$ onto $A$
induced by the $\nu_i$'s. Also, the group $\prod_{i\in I}^{\st{n}{*}}A_i$ will be assumed to be the $n$th nilpotent product of $A_i$'s. If $\psi_n$ is the natural homomorphism from $A$ onto $G_n$ induced by the identity map on each $A_i$,
$$ L=\prod_{i\in I}^*L_i\st {\nu }{\lra }A=\prod_{i\in I}^*A_i\st {\psi_n}{\lra}
G_n=\prod_{i\in I}^{\st {n}{*}}A_i\lra 1\ ,$$
then we denote by $H_n$ the kernel of $\psi_n$ and $P_n$ the inverse image of
$H_n$ in $L$ under $\nu$; $ H_n=ker\psi_n \ \ \& \ \ P_n=\nu^{-1}(H_n)\ .$
We also put $D_c=\prod_{\exists j,\mu_j\neq i}[M_i,L_{\mu_1},\ldots ,L_{\mu_c}]^L$ and
$ E_c=D_1\cap \ga_{c+1}(L).$
Now, by a similar proof of Lemma 3.1 of [3] we have \\
{\bf Lemma A}

$(i)\ H_n=\nu (J_n)$ , where $J_n=\ga_{n+1}(L)\cap [L_i]^*$, ($[L_i]^*$ is the cartesian subgroup of the free product $L=\prod_{i\in I}^{*}L_i$ (see [5]).)\\

$(ii)\ G_n\cong L/P_n$ , where $P_n=ker(\psi_n \circ \nu)=(\prod_{i\in I}
M_i^L)J_n \ .$\\

 By rewriting the proof of the main result of section $3$ of [3] and noticing that
$$ {\cal N}_cM(A_i)\cong \frac{M_i\cap \ga_{c+1}(L_i)}{[M_i,\ _cL_i]}\ \ ,$$
when $L_i$ is a free group and also paying attention to this fact that being free for $F_i$'s  is not necessary in the proof, we can state the following vital lemma.\\
{\bf Lemma B}

 By the above notation and assumption,\\
$(i)$ If $n<c$, then
$$ \frac {(\prod_{i}^{}(M_i\cap \ga_{c+1}(L_i))E_cJ_c}{\prod_{i}^{}[M_i,\ _cL_i]E_c[J_n,\
_cL]}\cong \prod_{i}^{}\!^{\times}\frac {M_i\cap
\ga_{c+1}(L_i)}{[M_i,\ _cL_i]}\oplus \frac {H_c}{[H_n,\ _cA]}\ .$$
$(ii)$ If $n\geq c$, then
$$ \frac {(\prod_{i}^{}(M_i\cap \ga_{c+1}(L_i))E_cJ_n}{\prod_{i}^{}[M_i,\ _cL_i]E_c[J_n,\
_cL]}\cong \prod_{i}^{}\!^{\times}\frac {M_i\cap
\ga_{c+1}(L_i)}{[M_i,\ _cL_i]}\oplus \frac {H_n}{[H_n,\ _cA]}\
.$$\\

 From now on, we assume that each $A_i$ is cyclic, therefore its
Baer-invariant is trivial i.e $M_i=1$ for all $i\in I$. Thus, we
have $P_n=J_n,\ D_c=E_c=1$.

 By considering natural epimorphisms $\theta_n:L\lra \frac {L}{[J_n,\ _cL]}$, we define
$$ \ol {S_n}=\theta_n(L)\ \ and\ \ \ol {P_n}=\theta_n(P_n)\ .$$
Clearly
$$ \ol{S_n}=\frac{L}{[J_n,\ _cL]}\ \ and \ \ \ol{P_n}\cong \frac {P_n}{P_n\cap [J_n,\ _cL]}=\frac {P_n}{[J_n,\ _cL]}\ .$$
Now, we are in a position to state and prove the main result of
this paper.\\
{\bf Theorem}

 Let $G_n=\prod_{i\in I}^{\st{n}{*}}A_i$ be the $n$th nilpotent product of cyclic groups $A_i$'s. Then by the above notation and assumption, we have\\
$(i)$ If $n\geq c$, then the following exact sequence is an
${\cal N}_c$-stem cover for $G_n$
$$ 1\lra \ol{P_n}\st {\subseteq }{\lra}\ol{S_n}\lra G_n\lra 1\ ,$$
i.e $\ol{P_n}\subseteq Z_c(\ol{S_n})\cap\ga_{c+1}(\ol{S_n})$ and $\ol{P_n}\cong {\cal N}_cM(G_n)$.
 In other words, $\ol{S_n}$ is an ${\cal N}_c$-covering group for $G_n$.\\
$(ii)$ If $n<c$, then $G_n$ is the only ${\cal N}_c$-covering group of itself when
${\cal N}_cM(G_n)=1$, otherwise, $G_n$ has no ${\cal N}_c$-covering group.\\
{\bf Proof} $(i)$ By Lemma $A$ $(ii)$, clearly
$$ G_n\cong \frac{L}{\ol{P_n}}\cong \frac {L/[J_n,\ _cL]}{P_n/[J_n,\ _cL]}\cong \frac {\ol{S_n}}
{\ol{P_n}}\ .$$
So the sequence $1\ra \ol{P_n}\ra \ol{S_n}\ra G_n\ra 1$ is exact. Since $P_n=J_n$, we have
$$[\ol{P_n},\ _c\ol{S_n}]=[\theta_n(P_n),\ _c\theta_n(L)]=\theta([P_n,\ _cL])=\theta([J_n,\ _cL])=1\ .$$
So $\ol{P_n}\leq Z_c(\ol{S_n})$. Also, since $n\geq c$,
$P_n=J_n\leq \ga_{n+1}(L)\leq \ga_{c+1}(L)$. Therefore
$$ \ol{P_n}=\theta_n(P_n)\leq \theta_n(\ga_{c+1}(L))=\ga_{c+1}(\theta_n(L))=\ga_{c+1}(\ol{S_n})\ .$$
Hence $\ol{P_n}\leq Z_c(\ol{S_n})\cap \ga_{c+1}(\ol{S_n})$. Moreover
$$\ol{P_n}\cong \frac{P_n}{[J_n,\ _cL]}=\frac{P_n\cap
\ga_{c+1}(L)}{[P_n,\ _cL]}\ \ \ \ (since\ P_n=J_n\ and\ n\geq c).$$
By Lemma 3.4 of [3], we have
$$\ol{P_n}\cong \frac{P_n\cap
\ga_{c+1}(L)}{[P_n,\ _cL]}=\frac {(\prod_{i}^{}(M_i\cap
\ga_{c+1}(L_i))E_cJ_n}{\prod_{i}^{}[M_i,\ _cL_i]E_c[J_n,\ _cL]}$$
$$ \cong  \prod_{i}^{}\!^{\times}\frac {M_i\cap
\ga_{c+1}(L_i)}{[M_i,\ _cL_i]}\oplus \frac {H_n}{[H_n,\ _cA]}\ \
\ (By\ Lemma\ B\ and\ D_c=E_c) $$ $$\cong \frac{H_n}{[H_n,\
_cA]}\ \ \ \ \ \ \ \ \ \ (since\ M_i=1\ for\ all\ i\in I)$$
$$\cong {\cal N}_cM(G_n)\ \ \ \ \ \ \ \ \ (By\ Theorem\ 4.1\ of\
[3]).$$

$(ii)$ If ${\cal N}_cM(G_n)=1$, then by definition, $G_n$ is the
only ${\cal N}_c$-covering group of $G_n$. Now, let ${\cal
N}_cM(G_n)\neq 1$. By definition of the $n$th nilpotent product
$$ G_n=\frac {A}{[A_i]^*\cap\ga_{n+1}(A)}\ .$$
Since $A_i$'s are cyclic, it is easy to see that
$\ga_{n+1}(A)=[A_i]^*\ga_{n+1}(A)$ (see [5]). Therefore $G_n$ is
a nilpotent group of class $n$. Hence, by Theorem 2.1 of [4],
$G_n$ has no ${\cal N}_c$-covering group. $\Box$

 Note that the above theorem is a generalization of a result of
 Haebich [1, Theorem 3.1].\\
 {\bf Corollary}

 For any natural number $n$, there exists a nilpotent group $G_n$
 of class $n$ such that its Baer-invariant with respect to the
 variety ${\cal N}_c$ is {\it not} trivial and has at least one
 ${\cal N}_c$-covering group for all $c\leq n$.\\
 {\bf Proof.} Let $\{A_i|i\in I\}$ be a family of infinite cyclic
 groups and $card(I)=m$. Then clearly the free product
 $\prod_{i\in I}{*}A_i$ is a free group of rank $m$ and the $n$th
 nilpotent product $\prod_{i\in I}^{\st{n}{*}}A_i=F/\ga_{n+1}(F)$
 is a free nilpotent group of class $n$. By definition, for any
 $c\leq n$ we have
$${\cal N}_cM(G_n)=\frac{\ga_{n+1}(F)}{[\ga_{n+1}(F),\ _cF]}=\frac{\ga_{n+1}(F)}{\ga_{n+c+1}(F)}\neq 1\ .$$
Hence ${\cal N}_cM(G_n)$ is a free abelian group of rank
$\sum_{i=n+1}^{n+c}\alpha_i$ where $\alpha_i$ is the number of
basic commutators of weight $i$ on $m$ letters. Now, by the
previous theorem part $(i)$ $G_n$ has at least one ${\cal
N}_c$-covering group. $\Box$

 Note that for any $c>n,\ G_n$ is a nilpotent group of class $n$ such that its Baer-invariant, ${\cal N}_cM(G_n)=\ga_{c+1}(F)/\ga_{n+c+1}(F)$, is a free abelian group of rank $\sum_{i=c+1}^{n+c}\alpha_i$ and hence by the main theorem part $(ii)$, has no ${\cal N}_c$-covering group.\\
{\bf Remark}

 The above corollary gives a partial answer to the question that arose in [4] about the existence of ${\cal N}_c$-covering group for nilpotent groups of class $n$, when $c\leq n$.

\end{document}